\documentclass{amsart}
\usepackage{color}
\usepackage{graphicx}
\usepackage{latexsym}
\usepackage{amsfonts}
\usepackage[all]{xy}
\usepackage{amssymb, mathrsfs, amsfonts, amsmath}
\usepackage{amsbsy}
\usepackage{amsfonts}
\newtheorem*{acknowledgement}{Acknowledgement}

\newtheorem{remark}{Remark}
\newtheorem{theorem}{Theorem}
\newtheorem{example}{Example}
\numberwithin{equation}{section}

\begin{document}
	\title[Invariant soltutions for the Einstein equation]{INVARIANT SOLTUTIONS FOR THE EINSTEIN FIELD EQUATION}
	\author{Marcelo Barboza$^{1}$}
	\address{$^{1}$ Universidade Federal de Goi\'as, IME, 131, 74001-970, Goi\^ania, GO, Brazil.}
	\email{salve.barboza@gmail.com $^{1}$}
	\author{Benedito Leandro$^{2}$}
	\address{$^{2}$ Universidade Federal de Goi\'as, Centro de Ci\^encias Exatas, Regional Jata\'i, BR 364, km 195, 3800, 75801-615, Jata\'i, GO, Brazil.}
	\email{bleandroneto@gmail.com$^{2}$}
	\author{Romildo Pina$^{3}$}
	\address{$^{3}$ Universidade Federal de Goi\'as, IME, 131, 74001-970, Goi\^ania, GO, Brazil.}
	\email{romildo@ufg.br$^{3}$}
	\thanks{$^{1}$ Supported by CAPES/Brazil Proc. No. 303774/2009-6  Ministry of Science and Technology Brazil.}
	\date{October 03, 2017}

	\begin{abstract}
		In this paper we provide a method capable of producing an infinite number of solutions for Einstein's equation on static spacetimes with perfect fluid as a matter field. All spacetimes of this type which are symmetric with respect to a given group of translations and whose spatial factor is conformally flat, are characterized. We use this method to give some exact solutions of the referred equation. 
	\end{abstract}
	
	\maketitle
	\section{Introduction and Main Results}
	
	The Einstein equation on a spacetime $\widehat{M}$ containing matter with stress-energy tensor $T$ is given by
	\begin{equation}\label{einstein_equation}
	\widehat{G}=8\pi T,
	\end{equation} where $\widehat{G}$ is the Einstein gravitational tensor of $\widehat{M}$ and $T$ is a symmetric $(0,2)$-tensor with divergence zero. Roughly speaking, this equation shows how matter determines the Ricci curvature of spacetime and the vanishing of the divergence of $T$ tells how Ricci curvature moves this matter. The manifold $\widehat{M}$ is said to be a vacuum spacetime in case $T=0$ (cf. \cite{LeandroPina,Romildo,Oneil}). 
	
	A semi-Riemannian $(n+1)$-manifold $(\widehat{M}^{n+1},\,\widehat{g})$ is called static if $\widehat{M}$ is the product manifold $M\times\mathbb{R}$ of a manifold $M^n$ with $\mathbb{R}$ and the metric $\widehat{g}$ has the form
	\begin{equation}\label{static_metric}
	\widehat{g}=x^*g-f(x)^2dt^2,
	\end{equation} where $x:\widehat{M}\to M$ and $t:\widehat{M}\to\mathbb{R}$ are the natural projections, $g$ is a semi-Riemannian metric on $M$ and $f\in C^{\infty}(M)$ is strictly positive.
	
	A perfect fluid stress-energy tensor may be recognized by the identity
	\begin{equation}\label{perfect_fluid_tensor}
	T=(\mu+\rho)\eta\otimes\eta+\rho g,
	\end{equation} where $\mu,\,\rho\in C^{\infty}(M)$ are the energy density and pressure, respectively, and $\eta$ is a $1$-form such that $g(\eta,\eta)=-1$ and whose associated vector field represents the flux of the fluid (cf. \cite{Obata,Oneil}). The definition of a perfect fluid does not provide, however, a way to construct a spacetime model of one. We consider Einstein's equation on $(\widehat{M},\,\widehat{g})$ with perfect fluid as a matter field. It implies that the Ricci tensor of $(\widehat{M},\,\widehat{g})$ has, at each point of $\widehat{M}$, at most two distinct eigenvalues with multiplicities $1$ and $n$ (cf. \cite{Obata}).
	
	%
	By means of (\ref{einstein_equation}), (\ref{static_metric}) and (\ref{perfect_fluid_tensor}) we are able to conclude that (see \cite{Obata})
	\begin{eqnarray}\label{equationP}
	f\left(Ric_{g}-\frac{R_{g}}{n}g\right)=\left(\nabla^{2}_{g}f-\frac{\Delta_{g}f}{n}g\right),
	\end{eqnarray} and
	\begin{eqnarray}\label{densidadeEpressão}
	\mu=\frac{R_{g}}{2},\,\,\,\, f\rho=\left(\frac{n-1}{n}\right)\left[\Delta_{g}f - \frac{(n-2)}{2(n-1)}R_{g}f\right],
	\end{eqnarray} where $Ric_{g}$, $R_{g}$, $\nabla^2_g$ and $\Delta_g$ are, respectively, the Ricci tensor, scalar curvature, Hessian and Laplacian with respect to $g$. Observe that both sides of (\ref{equationP}) are trace-free. Kobayashi and Obata \cite{Obata}, proved that a perfect fluid static spacetimes whose spatial factor is conformally flat and complete admits a special warped product structure (cf. Lemma B in \cite{kobayashi}). On a static vacuum spacetime, which corresponds to $\mu=\rho=0$, relations (\ref{equationP}) and (\ref{densidadeEpressão}) together give
	\begin{eqnarray}\label{equationPvacuum}
	fRic_{g}=\nabla^{2}_{g}f\quad\mbox{and}\quad\Delta_{g}f=0,
	\end{eqnarray}
	telling that $(\widehat{M}^{n+1},\,\hat{g})$ must be a Ricci-flat manifold.
	
	The Einstein equation represents a great step when it comes to comprehend the Universe and the exact solutions of this equation have played a crucial role in the development of cosmology. Many scientists contributed with their spacetime models, like Schwarzchild, Kerr, Curzon, de Sitter, Bach, Levi-Civita and Kinnersly (cf. \cite{Hawking,Kramer,Oneil}). Exact solutions of Einstein's equations thus model gravitating systems and leads to the mathematical and physical understanding of the general theory of relativity.
	
	In this paper we construct perfect fluid static spacetimes whose spatial factor is conformally flat and admits a group of symmetries consisting of translations (cf. \cite{Olver}). Leandro, Pina and Sousa have recently applied this method to obtain static vacuum spacetime models (cf. \cite{LeandroPina,Romildo}). Specifically, we consider a pseudo-Euclidean metric \[\delta=\sum^n_{i=1}\varepsilon_idx_i\otimes dx_i\] in Cartesian coordinates $x=(x_1,\ldots,x_n)$ of $\mathbb{R}^n$ where $n\geq3$, $\varepsilon_1,\ldots,\varepsilon_n=\pm1$ and $\varepsilon_i=1$ for at least one $i\in\{1,\ldots,n\}$. Take $\xi:\mathbb{R}^n\to\mathbb{R}$ by the expression \[\xi(x_1,\ldots,x_n)=\alpha_1x_1+\cdots+\alpha_nx_n,\] so defined for an arbitrary choice of $\alpha=(\alpha_1,\ldots,\alpha_n)\in\mathbb{R}^n\setminus\{0\}$ and every $(x_1,\ldots,x_n)\in\mathbb{R}^n$. We then look for smooth functions $\varphi,\,f:(a,b)\subset\mathbb{R}\to(0,\infty)$ such that the composites $\varphi=\varphi\circ\xi,\,f=f\circ\xi:M\subset\mathbb{R}^n\to (0,\infty)$, where $M=\{x\in\mathbb{R}^n:a<\xi(x)<b\}=\xi^{-1}(a,b)$ is open in $\mathbb{R}^n$, both satisfy (\ref{equationP}) and (\ref{densidadeEpressão}). Therefore, the manifold $\widehat{M}^{n+1}=M\times\mathbb{R}$ with metric tensor \[\widehat{g}=\frac{\delta}{\varphi(x)^2}-f(x)^2dt^2\] is a perfect fluid static spacetime, symmetric with respect to the additive group $G=\{x\in\mathbb{R}^n:\xi(x)=0\}$ of translations in $\mathbb{R}^n$ and whose spatial factor is conformally flat. What has been said above is summed up in the next result.
	
	\begin{theorem}\label{theorem2}
		With $(\mathbb{R}^n,\,\delta)$ and $\varphi=\varphi\circ\xi,\,f=f\circ\xi$ as above, the manifold $\widehat{M}=M\times\mathbb{R}$, furnished with the metric tensor \[\widehat{g}=\frac{\delta}{\varphi(x)^{2}}-f(x)^2dt^2,\] is a perfect fluid static spacetime if, and only if,
		\begin{equation}\label{EDO}
		(n-2)f\varphi''- f''\varphi-2\varphi'f'=0
		\end{equation}
		where the energy density and pressure functions are given by
		\begin{eqnarray}\label{DenEdo}
		\mu(\xi)&=&\|\alpha\|^{2}\left(n-1\right)\left[\varphi\varphi''-\frac{n}{2}(\varphi')^{2}\right]
		\end{eqnarray}
		%
		and
		\begin{eqnarray}\label{preEdo}
		\rho(\xi)&=&\|\alpha\|^{2}\left(\frac{n-1}{n}\right)\left[\varphi^{2}\frac{f''}{f}-(n-2)\left(\frac{f'}{f}\varphi'\varphi + \varphi\varphi''-\frac{n}{2}(\varphi')^{2}\right)\right],
		\end{eqnarray} 
		respectively, where $\|\alpha\|^{2}=\sum^n_{k=1}\varepsilon_k\alpha^2_k$ equals $-1$, $0$ or $1$ if $\alpha$ is a timelike, lightlike or spacelike vector, respectively.
	\end{theorem}
	
	The Riemannian solutions for Theorem \ref{theorem2} are not complete (cf. \cite{kobayashi,Obata}). However, we prove a semi-Riemannian example which is geodesically complete (see Example \ref{exampleComplet}).
	
	\begin{remark}\label{riccatti_equation}
		Upon dividing both sides of (\ref{EDO}) by $\varphi f$ we get
		\begin{eqnarray}\label{1.3.1}
		(n-2)\frac{\varphi^{\prime\prime}}{\varphi}=\frac{f^{\prime\prime}}{f}+2\frac{\varphi^{\prime}}{\varphi}\frac{f^{\prime}}{f},
		\end{eqnarray} or even,
		\begin{eqnarray}\label{riccatti}
		y^{\prime}=(n-2)(x^{\prime}+x^2)-2xy-y^2
		\end{eqnarray}
		where $x=\frac{\varphi^{\prime}}{\varphi},\,y=\frac{f^{\prime}}{f}$. The above equation being Riccati in $y$ for a known $x$ (cf. \cite{Ince}), its solutions are written in the form $y=y_0+u$, where $y_0$ is a particular solution of (\ref{riccatti}) and $u$ satisfies \[\frac{d}{d\xi}\left(\frac{1}{u}\right)-2(y_0+x)\frac{1}{u}=1.\]
	\end{remark}
	
	Theorem \ref{theorem2} says that there is a perfect fluid static spacetime of distinguished geometric nature corresponding to each direction $\alpha\in\mathbb{R}^n\setminus\{0\}$ and smooth $\varphi:(a,b)\subset\mathbb{R}\to(0,\infty)$. The authors of \cite{BarbosaPinaTenenblat,LeandroPina,LeandroTenenblat,LeandroPinaSousa,Romildo} verified a similar phenomenon under the assumption that $\|\alpha\|^2=0$.
	
	For the remaining of this section we will apply Theorem \ref{theorem2} on the construction of $5$ different spacetime models by either specifying the conformal factor $\varphi$ or the energy density $\mu$. It is worth noticing that $\mu>|\rho|$ for all known forms of matter (cf. \cite{Oneil}). Therefore, some of our examples do not have an apparent physical meaning. It is a remarkable fact that some manifolds do not admit geodesically complete Riemannian metrics of non negative scalar curvature (cf. \cite{gromov1983positive}). Since $\mu$ equals half the scalar curvature of the spacetime's spatial factor, finding a realistic model is a real challenge.
	
	\begin{example}
		If we take $\mu=(\varphi^{\prime})^{2}>0$ then, according to (\ref{DenEdo}), we have
		\begin{eqnarray*}
			\|\alpha\|^{2}(n-1)\left[\varphi\varphi^{\prime\prime}-\frac{n}{2}(\varphi^{\prime})^{2}\right]=(\varphi^{\prime})^{2}
		\end{eqnarray*}
		and hence,
		\begin{eqnarray*}
			\frac{\varphi^{\prime\prime}}{\varphi^{\prime}}=\left[\frac{1+\|\alpha\|^{2}(n-1)\frac{n}{2}}{\|\alpha\|^{2}(n-1)}\right]\frac{\varphi^{\prime}}{\varphi}.
		\end{eqnarray*}
		Upon integrating twice with respect to $\xi$ we get
		\[\varphi(\xi)=\left[\eta(\tilde{\kappa}\xi+\bar{\kappa})\right]^{\frac{1}{\eta}}\]
		where $\tilde{\kappa}\neq0$, $\bar{\kappa}$ and $\eta=\left(1-\frac{2+\|\alpha\|^{2}n(n-1)}{2\|\alpha\|^{2}(n-1)}\right)$ are constants. By setting \[n=3,\quad\|\alpha\|^{2}=1,\quad\tilde{\kappa}=1\quad\mbox{and}\quad\bar{\kappa}=0,\] we find \[\varphi(\xi)=-\xi^{-1}\] and (\ref{EDO}) then becomes
		\[f^{\prime\prime}-\frac{2}{\xi}f^{\prime}-\frac{2}{\xi^{2}}f=0,\] of whose \[f(\xi) = \xi^{\frac{3-\sqrt{17}}{2}}\] is a solution. A straightforward computation gives \[\mu(\xi)=\frac{1}{\xi^4}>\rho(\xi)=\left(4-\sqrt{17}\right)\frac{1}{\xi^{4}}\] and the manifold \[\widehat{M}^{n+1}=M\times\mathbb{R},\quad M=\{x\in\mathbb{R}^n:\xi(x)>0\},\] with metric
		\begin{eqnarray*}
			\hat{g}=(\alpha_{1}x_{1}+\alpha_{2}x_{2}+\alpha_{3}x_{3})^{2}(dx_{1}^{2}+dx_{2}^{2}+dx_{3}^{2})-\frac{(\alpha_{1}x_{1}+\alpha_{2}x_{2}+\alpha_{3}x_{3})^{3}}{(\alpha_{1}x_{1}+\alpha_{2}x_{2}+\alpha_{3}x_{3})^{\sqrt{17}}}dt^{2}
		\end{eqnarray*}
		is a perfect fluid static spacetime with energy and pressure, respectively, given by
		\begin{eqnarray*}
			\mu=\frac{1}{(\alpha_{1}x_{1}+\alpha_{2}x_{2}+\alpha_{3}x_{3})^{4}}\quad\mbox{and}\quad\rho=\frac{4-\sqrt{17}}{(\alpha_{1}x_{1}+\alpha_{2}x_{2}+\alpha_{3}x_{3})^{4}}.
		\end{eqnarray*}
	\end{example}
	
	\vspace{8pt}
	
	\begin{example}\label{exampleComplet}
		Taking $\varphi(\xi)=e^{\xi}$ on equation (\ref{1.3.1}) produces \[f^{\prime\prime}+2f^{\prime}-(n-2)f=0,\] which is a second order linear equation with constant coefficients and, as such, it's solutions must be linear combinations of the functions \[f^{+}(\xi)=e^{(1+\sqrt{n-1})\xi}\quad\mbox{and}\quad f^{-}(\xi)=e^{(1-\sqrt{n-1})\xi}.\]
		Therefore, the manifold \[\widehat{M}=\mathbb{R}^n\times\mathbb{R}\] with metric tensor \[g=\frac{\delta}{e^{2\xi}}-e^{2(1+\sqrt{n-1})\xi}dt^2\] furnishes an example of a perfect fluid static spacetime with energy and pressure given by \[\mu(\xi)=-\frac{\|\alpha\|^{2}(n-1)(n-2)}{2}e^{2\xi}\] and \[\rho(\xi)=\|\alpha\|^{2}\left(\frac{n-1}{n}\right)\left\{\frac{(n-2)^{2}}{2}-(1+\sqrt{n-1})(n-3-\sqrt{n-1})\right\}e^{2\xi},\] respectively. In case that $\alpha$ is a lightlike vector field, this is a geodesically complete spacetime (cf. \cite{LeandroPina}).
	\end{example}
	
	\vspace{8pt}
	
	\begin{example}
		The previous example is enlightening when it shows us that there are certain cases for which $f(\xi)$ may be writen in terms of $\varphi(\xi)$, as in
		\[f^{+}(\xi)=e^{(1+\sqrt{n-1})\xi}=\varphi(\xi)^{(1+\sqrt{n-1})}.\]
		Henceforth, we keep our duty of looking for solutions of (\ref{1.3.1}), but this time with the further assumption that $f$ is a function of $\varphi$, such as $f(\xi)=\frac{1}{\varphi(\xi)}$. In this case, equation (\ref{EDO}) becomes \[(n-1)\frac{\varphi^{\prime\prime}}{\varphi}=0\iff\varphi^{\prime\prime}=0,\] forcing $\varphi$ to be a linear function of $\xi$, i.e., $\varphi(\xi)=a\xi+b$ where $a,b\in\mathbb{R}$ are constants, $a\neq0$ and $\varphi(\xi)$ is defined on the open half-space $M=\{x\in\mathbb{R}^n:a\xi(x)+b>0\}$ so that it is a positive function. The perfect fluid static spacetime thus obtained is described as the manifold \[\widehat{M}^{n+1}=M\times\mathbb{R}\] alongside metric tensor \[\widehat{g}=\frac{\sum_k\varepsilon_kdx^2_k-dt^2}{(a\sum_k\alpha_kx_x+b)^2}\] and constant energy and pressure functions given by \[\mu=-\frac{n(n-1)\|\alpha\|^{2}}{2}\quad\mbox{and}\quad\rho=\frac{n(n-1)\|\alpha\|^{2}}{2},\] respectively. The manifold $(\widehat{M}^{n+1},\widehat{g})$ is, therefore, globally conformal to a pseudo-Euclidean space whose dimension and signature have both been increased
		by $1$ for those starting with $(M,g)$. Special instances of the above setting include
		\[\varepsilon_1=\cdots=\varepsilon_n=1,\,\,\,\alpha_1=\cdots=\alpha_{n-1}=0,\,\,\,\alpha_n=1,\,\,\,a=1\quad\mbox{and}\quad
		b=0,\] corresponding to the Lorentzian manifold
		\[\varphi(\xi)=\xi=x_n,\quad f(\xi)=\frac{1}{x_n},\]
		\[M=\{(x_1,\ldots,x_n)\in\mathbb{R}^n:x_n>0\}=\mathbb{H}^n,\quad\widehat{M}^{n+1}=\mathbb{H}^n\times\mathbb{R}\]
		with metric \[\widehat{g}=\frac{\sum_kdx^2_k-dt^2}{x_n^2},\] whose spatial factor is the $n$-dimensional hyperbolic space.  
	\end{example}
	
	\vspace{8pt}
	
	Although we can choose $\varphi$ at (\ref{EDO}), it seems that finding an explicit formula for $f$ is a difficult task.
	
	\vspace{8pt}
	
	\begin{example}
		If we take $x(\xi)=-\tan(\xi)$ at (\ref{riccatti}), the referred equation shows itself up like \[y^{\prime}=-(n-2)+\tan(\xi)^2-(\tan(\xi)-y)^2,\] and even, like \[z^{\prime}=-(z^2+n-1),\] on the variable $z(\xi)=y(\xi)-\tan(\xi)$. Since the general solution of the above equation is
		\[z(\xi)=-\sqrt{n-1}\tan(a+\xi\sqrt{n-1})\\[0.3cm]\] where $a\in\mathbb{R}$ is some constant, we get from
		\[\frac{f^{\prime}(\xi)}{f(\xi)}=y(\xi)=z(\xi)+\tan(\xi)\] that
		\[f(\xi)=b\frac{\cos(a-\xi\sqrt{n-1})}{\cos(\xi)}\] for some constant $b>0$. Also, it is readily seen from \[\frac{\varphi^{\prime}(\xi)}{\varphi(\xi)}=x(\xi)=-\tan(\xi)\] that \[\varphi(\xi)=c\cdot\cos(\xi)\] for some constant $c\in\mathbb{R}\setminus\{0\}$. By taking \[n=3,\quad a= 0,\quad b=1\quad\mbox{and}\quad c=1,\] the perfect fluid static spacetime $(\widehat{M}^{n+1},\widehat{g})$ that we have got might then be pictured as the manifold
		\[\varphi(\xi)=\cos(\xi),\quad f(\xi)=\frac{\cos(\xi\sqrt{2})}{\cos(\xi)},\]
		\[\widehat{M}=M\times\mathbb{R},\quad M=\left\{x\in\mathbb{R}^n:-\frac{\sqrt{2}\pi}{4}<\xi(x)<\frac{\sqrt{2}\pi}{4}\right\},\] along with metric
		\[\widehat{g}=\frac{\varepsilon_1dx^2_1+\varepsilon_2dx^2_2+\varepsilon_3dx^2_3-\cos^2(\sqrt{2}(\alpha_{1}x_{1}+\alpha_{2}x_{2}+\alpha_{3}x_{3}))dt^2}{\cos^2(\alpha_{1}x_{1}+\alpha_{2}x_{2}+\alpha_{3}x_{3})}\] and functions
		\[\begin{array}{rcl}
		\mu(\xi)&=&\|\alpha\|^{2}\left[\cos^{2}(\xi)-3\right]\quad\mbox{and} \\ [0.3cm]
		\rho(\xi)&=&-\frac{\|\alpha\|^{2}}{3}[2\sqrt{2}\tan(\sqrt{2}\xi)\cos(\xi)(2\sin(\xi)-1)\nonumber\\
		&+&7\cos^{2}(\xi)+2\sin(\xi)-7] \\ [0.3cm]
		\end{array}\]
		serving as the energy density and pressure of the fluid, respectively.
	\end{example}
	
	\vspace{8pt}
	
	\begin{example}
		Back at the Euclidean scenario, let's choose $\mu(\xi)=\frac{n-1}{2}\|\alpha\|^2\varphi(\xi)^2>0$. From (\ref{DenEdo}) we have \[\left(\frac{n-1}{2}\right)\|\alpha\|^2\left\{2\varphi\varphi^{\prime\prime}-n(\varphi^{\prime})^2\right\}=\left(\frac{n-1}{2}\right)\|\alpha\|^2\varphi^2\] and upon dividing both sides of the above equation by $\frac{n-1}{2}\|\alpha\|^2\varphi^2$ we conclude that \[2{\frac{\varphi^{\prime\prime}}{\varphi}}-n\left(\frac{\varphi^{\prime}}{\varphi}\right)^2=1,\] the last being equivalent to \[2\frac{d}{d\xi}\left(\frac{\varphi^{\prime}}{\varphi}\right)-(n-2)\left(\frac{\varphi^{\prime}}{\varphi}\right)^2=1.\] The general solution of the above equation is \[\varphi(\xi)=b\sec^{\frac{2}{n-2}}\left(a+\xi\frac{\sqrt{n-2}}{2}\right),\] where $a,b\in\mathbb{R}$, $b>0$. By setting \[n=3,\quad a=0,\quad\mbox{and}\quad b=1,\] we get \[\varphi(\xi)=\sec^2\left(\frac{\xi}{2}\right),\quad \mu(\xi)=\sec^4(\frac{\xi}{2}),\] and (\ref{EDO}) then becomes
		\begin{eqnarray}\label{ult2}
		f^{\prime\prime}+2\tan\left(\frac{\xi}2\right)f^{\prime}-\frac{1}{2}\left\{1+3\tan^2\left(\frac{\xi}2\right)\right\}f=0.
		\end{eqnarray}
		The phase portrait of (\ref{ult2}) is pictured in Figure \ref{PhasePortrait}.
		\begin{figure}[h!]
			\centering
			\includegraphics[width=0.4\textwidth]{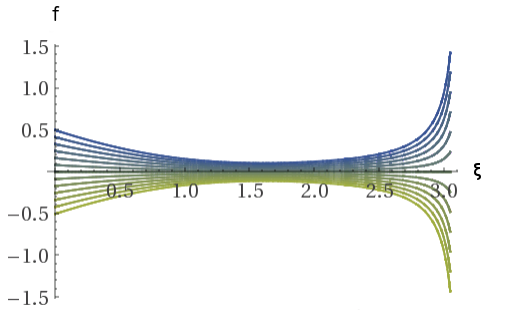}
			\caption{Sampling $f(0)$ and $f^{\prime}(0)$.}
			\label{PhasePortrait}
		\end{figure}
	\end{example}

	\section{Proof of the Main Result}
	
	In this section we prove Theorem \ref{theorem2}. We closely follow the lines in (cf. \cite{BarbosaPinaTenenblat,LeandroPina,LeandroTenenblat,LeandroPinaSousa,Romildo}).
	
	\
	
	\noindent {\bf Proof of Theorem \ref{theorem2}:}
	Let $\xi:\mathbb{R}^n\to\mathbb{R}$ be defined by the rule \[\xi(x_1,\ldots,x_n)=\alpha_1x_1+\cdots+\alpha_nx_n\] for every point $(x_1,\ldots,x_n)\in\mathbb{R}^n$ and some arbitrarily chosen direction $\alpha=(\alpha_1,\ldots,\alpha_n)\in\mathbb{R}^n\setminus\{0\}$. Take smooth functions $\varphi,\,f:(a,b)\subset\mathbb{R}\to(0,\infty)$ so we can make composites $\varphi=\varphi\circ\xi,\,f=f\circ\xi:M\subset\mathbb{R}^n\to(0,\infty)$, where $M=\{x\in\mathbb{R}^n:a<\xi(x)<b\}=\xi^{-1}(a,b)$ is open in $\mathbb{R}^n$. Also, consider a pseudo-Euclidean metric \[\delta=\sum_i\varepsilon_idx_i\otimes dx_i\] in coordinates $x=(x_1,\ldots,x_n)$ of $M\subset\mathbb{R}^n$, where $\varepsilon_1,\ldots,\varepsilon_n=\pm1$ and $\varepsilon_i=1$ for at least one $i$. The manifold $\widehat{M}=M\times\mathbb{R}$ with metric tensor \[\widehat{g}=g-f^2dt^2,\] where $g=\varphi^{-2}\delta$ is in the conformal class of $\delta$, is a solution of the Einstein equation with perfect fluid as a matter field if, and only if, $\varphi,\,f$ satisfy (\ref{EDO}) and $\mu,\,\rho$ are as in (\ref{densidadeEpressão}) (cf. \cite{Obata}). Since we have (see \cite{Besse})
	\begin{eqnarray}\label{conformal_ricci_tensor}
	Ric_{g}=\frac{1}{\varphi^2}\left\{(n-2)\varphi\nabla^2\varphi+[\varphi\Delta\varphi-(n-1)\|\nabla\varphi\|^2]\delta\right\},
	\end{eqnarray}
	where $\nabla\varphi$, $\nabla^2\varphi$ and $\Delta\varphi$ stand for the gradient, Hessian and Laplacian of $\varphi$ with respect to $\delta$, respectively. The
	scalar curvature
	$R_{g}=\displaystyle\sum_k\varphi^2\varepsilon_k(Ric_{g})_{kk}$ of
	$g$ satisfies
	\begin{eqnarray}\label{conformal_scalar_curvature}
	R_{g}=(n-1)\left\{2\varphi\Delta\varphi-n\|\nabla\varphi\|^2\right\},
	\end{eqnarray}
	and hence
	\begin{eqnarray}\label{conformal_trace_free_ricci_tensor}
	Ric_{g}-\frac{R_{g}}{n}g=\frac{n-2}{\varphi}\left(\nabla^2\varphi-\frac{\Delta\varphi}{n}\delta\right).
	\end{eqnarray}
	Being $\delta$ a flat metric tensor, it produces simple formulas like
	\[\begin{array}{rcl}
	(\nabla^2\varphi)_{ij} &=& \varphi_{x_ix_j} \\ [0.3cm]
	&=&\varphi^{\prime\prime}\xi_{x_i}\xi_{x_j}+\varphi^{\prime}\xi_{x_ix_j} \\ [0.3cm]
	&=&\alpha_i\alpha_j\varphi^{\prime\prime}
	\end{array}\] and
	\[\begin{array}{rcl}
	\Delta\varphi &=& \sum_k\varepsilon_k\varphi_{x_kx_k} \\ [0.3cm]
	&=& \sum_k\varepsilon_k\alpha^2_k\varphi^{\prime\prime} \\ [0.3cm]
	&=&\|\alpha\|^2\varphi^{\prime\prime},
	\end{array}\] from which we detect the tensor components \[\left(Ric_g-\frac{R_g}ng\right)_{ij}=\frac{(n-2)\varphi^{\prime\prime}}{\varphi}\left(\alpha_i\alpha_j-\delta_{ij}\varepsilon_i\frac{\|\alpha\|^2}{n}\right).\] In order to compute the Hessian $\nabla^2_gf$ of $f$ relatively to $g$ we evoke the expression
	\begin{eqnarray}\label{conformal_hessian_of_f}
	(\nabla^2_{g}f)_{ij}=f_{x_ix_j}-\sum_k\Gamma^k_{ij}f_{x_k},
	\end{eqnarray}
	where the functions
	\[\begin{array}{rcl}
	\Gamma^k_{ij}&=&\frac{1}{2}\sum_k\left\{(g_{jl})_{x_i}+(g_{li})_{x_j}-(g_{ij})_{x_l}\right\}g^{lk}
	\\ [0.3cm]
	&=&
	-\frac{\varepsilon_k}{\varphi}\left\{\delta_{jk}\varepsilon_j\varphi_{x_i}+\delta_{ki}\varepsilon_k\varphi_{x_j}-\delta_{ij}\varepsilon_i\varphi_{x_k}\right\},
	\\ [0.3cm]
	\end{array}\] that is, $${\Gamma}^{k}_{ij}=0,\quad{\Gamma}^{i}_{ij}=-\frac{\varphi_{x_{j}}}{\varphi},
	\quad{\Gamma}^{k}_{ii}=\varepsilon_{i}\varepsilon_{k}\frac{\varphi_{x_{k}}}{\varphi},\quad{\Gamma}^{i}_{ii}=-\frac{\varphi_{x_{i}}}{\varphi},$$
	are the Christoffel symbols of $g$, $i,j,k=1,\ldots,n$. Therefore,
	\[\begin{array}{rcl}\label{conformal_hessian_of_f}
	(\nabla^2_{g}f)_{ij} &=& f_{x_ix_j}+\frac{\varphi_{x_i}f_{x_j}+\varphi_{x_j}f_{x_i}}{\varphi}-\delta_{ij}\varepsilon_i\sum_k\varepsilon_k\frac{\varphi_{x_k}f_{x_k}}{\varphi} \\ [0.3cm]
	&=& \alpha_i\alpha_jf^{\prime\prime}+(2\alpha_i\alpha_j-\delta_{ij}\varepsilon_i\|\alpha\|^2)\frac{\varphi^{\prime}f^{\prime}}{\varphi}
	\end{array}\]
	and the Laplacian $\Delta_{g}f=\sum_k\varphi^2\varepsilon_k(\nabla^2_{g}f)_{kk}$
	of $f$ with respect to $g$ is \[\Delta_{g}f=\varphi^2\|\alpha\|^2\left\{f^{\prime\prime}-(n-2)\frac{\varphi^{\prime}f^{\prime}}{\varphi}\right\}\] thus giving \[\left(\nabla^2_gf-\frac{\Delta_gf}{n}g\right)_{ij}=\left(f^{\prime\prime}+2\frac{\varphi^{\prime}f^{\prime}}{\varphi}\right)\left(\alpha_i\alpha_j-\delta_{ij}\varepsilon_i\frac{\|\alpha\|^2}{n}\right).\]
	We claim that $\alpha_i\alpha_j-\delta_{ij}\varepsilon_i\frac{\|\alpha\|^2}{n}\neq0$ for some choice of (possibly equal) indexes $i,j$. As a matter of fact, if $\alpha_i\alpha_j=0$ for all $i\neq j$, then $\alpha=\alpha_{i_0}e_{i_0}$
	for some $i_0\in\{1,\ldots,n\}$
	and \[\alpha^2_{i_0}-\varepsilon_{i_0}\frac{\|\alpha\|^2_g}{n}=\alpha^2_{i_0}-\varepsilon_{i_0}\frac{\varepsilon_{i_0}\alpha^2_{i_0}}{n}=\alpha^2_{i_0}\frac{n-1}{n}\neq0,\] since $\alpha\neq0$. If $\alpha_i\alpha_j\neq0$ for some $i\neq j$,
	then \[\alpha_i\alpha_j-\delta_{ij}\varepsilon_i\frac{\|\alpha\|^2}{n}=\alpha_i\alpha_j\neq0.\]
	Therefore, equation (\ref{EDO}) translates itself into
	\[(n-2)f\frac{\varphi^{\prime\prime}}{\varphi}=f^{\prime\prime}+2\frac{\varphi^{\prime}}{\varphi}f^{\prime}.\] As for the remaining identities, we have \vspace{12pt}
	\[\begin{array}{rcl}
	\mu&=&\frac{n-1}{2}\left(2\varphi\Delta\varphi-n\|\nabla\varphi\|^{2}\right) \\ [0.3cm]
	&=& (n-1)\|\alpha\|^2\left(\varphi\varphi^{\prime\prime}-\frac{n}{2}\left(\varphi^{\prime}\right)^2\right), \\ [1cm]
	\rho &=& \frac{n-1}{n}\left(\frac{\Delta_gf}{f}-\frac{n-2}{2(n-1)}R_g\right) \\ [0.3cm]
	&=& \frac{n-1}{n}\left\{\|\alpha\|^2\varphi^{2}\left[f^{\prime\prime}-(n-2)\frac{\varphi^{\prime}f^{\prime}}{\varphi}\right]
	-\frac{(n-2)\|\alpha\|^2}{2}\left[2\varphi\varphi^{\prime\prime}-n(\varphi^{\prime})^{2}\right]f\right\} \\ [0.3cm]
	&=&\frac{(n-1)\|\alpha\|^2}{n}\left[\varphi^{2}\frac{f''}{f}-(n-2)\left(\frac{f'}{f}\varphi'\varphi + \varphi\varphi''-\frac{n}{2}(\varphi')^{2}\right)\right],
	\end{array}\] which concludes this demonstration.
	\hfill $\Box$
	
	\

\begin{thebibliography}{BB}
		
		\bibitem{BarbosaPinaTenenblat} Barbosa, E., Pina, R., Tenenblat, K.: On gradient Ricci solitons conformal to a pseudo-Euclidian space. Israel Journal of Mathematics 200 (2014), 213-224.
		
		\bibitem{Besse} Besse, Arthur L.: Einstein Manifolds. Springer Berlin Heidelberg, 1987.
		
		\bibitem{gromov1983positive} Gromov, M., Lawson, H.B.: Positive scalar curvature and the Dirac operator on complete Riemannian manifolds. Publications math{\'e}matiques de l'IH{\'E}S, 58, Springer, (1983) 83-196.
		
		\bibitem{Ince} Hale, J.K.: Ordinary DiffentiaI Equations. Krieger Publishing Company, Malabar, Florida, 2° ed., 1980.
		
		\bibitem{Hawking} Hawking S.W., Ellis, G.F.R.: The Large Scale Structure of Space-Time. Cambridge Univ. Press, 1973.
		
		\bibitem{LeandroPina} Leandro, B., Pina, R.: Invariant solutions for the static vacuum equation. Journal of Mathematical Physics 58, 072502 (2017).
		
		\bibitem{LeandroTenenblat} Leandro, B., Tenenblat, K.: On gradient Yamabe solitons conformal to a pseudo-Euclidian space. dx.doi.org/10.1016/j.geomphys.2017.07.020.
		
		\bibitem{LeandroPinaSousa} Leandro, B., Pina, R., Sousa Lemes, M.: On the structure of Einstein warped product semi-Riemannian manifolds. arXiv:1708.04720v1 [math.DG] 15 Aug 2017.
		
		\bibitem{Romildo} Lemes de Sousa, M., Pina, R.: A family of warped product semi-Riemannian Einstein
		metrics. Differential Geometry and its Applications, 50, (2017) 105-115.
		
		\bibitem{Kramer} Kramer, D.: Exact solutions of Einstein's field equations. Cambridge: Cambridge University Press, 1980.
		
		\bibitem{kobayashi} Kobayashi, O.: A diferential equation arising from scalar curvature function. J. Math. Soc. Japan. Vol. 34, No. 4, 1982.
		
		\bibitem{Obata} Kobayashi, O., Obata M.:  Conformally-flatness and static space-time, in: Manifolds and Lie Groups, in: Progress in Mathematics, 14, Birkhäuser, (1981) 197-206.
		
		
		\bibitem{Olver} Olver, P. J.: Applications of Lie Groups to Differential Equations. Second Edition, GTM 107 (2000).
		
		\bibitem{Oneil} O'Neill, B.: Semi-Riemannian geometry with applications to relativity. Academic Press (1983).
		
		
	\end{thebibliography}
\end{document}